\documentclass[reqno,12pt]{amsart}
\newtheorem{Theorem}{Theorem}[section]
\newtheorem{Lemma}[Theorem]{Lemma}
\newtheorem{Proposition}[Theorem]{Proposition}

\newtheorem{Remark}[Theorem]{Remark}

\newtheorem{Definition}[Theorem]{Definition}

\renewcommand{\P}{{\mathbb P}} 
\renewcommand{\O}{{\mathcal O}} 

\newcommand{\LL}{{\mathcal L}} 
\newcommand{\MM}{{\mathcal M}} 
\newcommand{\GG}{{\mathcal G}} 
\newcommand{\EE}{{\mathcal E}} 
\newcommand{\BB}{{\mathcal B}} 
\newcommand{\Z}{{\mathbf  Z}} 
\newcommand{\T}{{\mathbb  T}} 
\newcommand{\GL}{GL} 
\newcommand{\tN}{{\widetilde N}}
\newcommand{\tM}{{\widetilde M}}

\newcommand{\uu}{{\underline u}}
\newcommand{\ux}{{\underline x}}
\newcommand{\up}{{\underline p}}
\newcommand{\uQ}{{\underline Q}}

\newcommand{\Sym}{\text{Sym}}
\newcommand{\Pic}{\text{Pic}}
\newcommand{\ch}{\text{ch}}
\newcommand{\td}{\text{td}}

\newcommand{\chk}{{$\looparrowleft$}}
\newcommand{\complex}{\mathbf C}
\newcommand{\proj}{\text{Proj}\;}
\newcommand{\rank}{\text{rank\,}}
\newcommand{\AS}{\mathcal A}
\newcommand{\Hilb}{\text{Hilb}}

\newcommand{\vspan}{\text{span}}
\newcommand{\image}{\text{image}}
\newcommand{\ra}{\rightarrow}

\newcommand{\lra}{\longrightarrow}

\renewcommand{\phi}{\varphi}

\newcommand{\demo}{\noindent {\sc Proof.}\;}

\usepackage{amsthm,amsmath,amsfonts,amscd,amssymb,latexsym,exscale,verbatim}

\begin{document}
\title{Apolar schemes of algebraic forms}
\author{Jaydeep Chipalkatti}
\maketitle

\parbox{12cm}{ \small 
{\sc Abstract.} 
This is a note on the classical Waring's problem for algebraic forms. 
Fix integers $(n,d,r,s)$, and let $\Lambda$ be a general $r$-dimensional 
subspace of degree $d$ homogeneous polynomials in $n+1$ variables. Let 
$\AS$ denote the variety of $s$-sided polar polyhedra of $\Lambda$. 
We carry out a case-by-case study of the structure of $\AS$ for several 
specific values of $(n,d,r,s)$. In the first batch of examples, $\AS$ is 
shown to be a rational variety. In the second batch, $\AS$ is a 
finite set of which we calculate the cardinality.}

\medskip 

\parbox{12cm}{\small
}

\bigskip 

\section{Introduction} 
We begin with a classical example to illustrate the theme of this paper. 
Let $F_1, F_2$ be general quadratic forms in variables $x_0, \dots, x_n$, 
with coefficients in $\complex$. It is then possible to diagonalize 
the $F_i$ simultaneously (see \cite[Ch.~22]{JoeH}), i.e., one can 
find linear forms $L_1, \dots, L_{n+1}$ such that 
\[ F_i = c_{i1}L_1^2 + \dots + c_{i(n+1)} L_{n+1}^2, \] 
for $i=1,2$, and some constants $c_{ij} \in \complex$. 
Moreover, up to rescaling there is a unique choice for 
the set $\{L_1, \dots, L_{n+1}\}$. This result naturally leads to 
similar questions about forms of higher degree, where much less 
is known in general. 

Now assume that $F_1,\dots, F_r$ are forms of degree  
$d$ in $x_0,\dots,x_n$. Let $Z = \{L_1,\dots, L_s\}$ be a  
collection of linear forms in the $x_i$, such that it is 
possible to write   
\[ F_i = c_{i1}L_1^d + \dots + c_{is} L_s^d, \qquad 1 \le i \le r; \]  
for some constants $c_{ij} \in \complex$. In nineteenth century 
terminology (introduced by Reye), $Z$ is then called a  
polar $s$-hedron (\emph{polar $s$-seit}) of the $\{F_i\}$. It corresponds to a 
collection of hyperplanes in $\P^n$ which stands in some geometric 
relation to the system of hypersurfaces defined by the $F_i$. The 
precise nature of this relation is very sensitive to the values 
$(n,d,r,s)$, but in any event it is invariant under the 
automorphisms of $\P^n$. 

For instance, in the example above, let $\Pi_i$ be the hyperplane 
defined by $L_i=0$. Then the $n+1$ points 
\[ P_k = \Pi_1 \cap \dots \cap \widehat{\Pi}_k \dots \cap \Pi_{n+1}, 
\quad (k=1,\dots,n+1) \] 
are exactly the vertices of the singular quadrics belonging to the pencil 
$\{ F_1 + \lambda F_2 =0 \}_{\lambda \in \P^1}$. 

\subsection{A Summary of Results} 
Fix degree $d$ forms $\{F_1, \dots, F_r \}$ as above. Then the polar 
$s$-hedra of this collection move in an algebraic family, denoted 
by $\AS$. (See Definition \ref{defn.polyhedron} \emph{et seq.~}for the precise 
statement.) In this note we deduce results about the birational structure 
of $\AS$ for several specific quadruples $(n,d,r,s)$, in each case 
assuming that the $F_i$ are chosen generally. 
A parameter count shows that the dimension of the variety $\AS$ is 
`expected' to be $s(n+r) -r \binom{n+d}{d}$ (more on this in 
\S \ref{section.prelim} below). 
For the quadruples 
\[ (2,4,2,8), \, (2,3,4,7), \, (3,2,6,7) \, (2,3,7,8), \] 
it is shown here that $\AS$ is a rational variety of expected dimension. 
For the cases 
\[ (2,3,8,8), \, (3,2,7,7), \, (2,4,3,9) \, (3,3,2,8), \, (2,3,3,6), \] 
the variety $\AS$ is expected to be (and is) a finite set of 
points; in each case we determine its cardinality. The 
calculation for $(2,3,3,6)$ was done by Franz London over a 
century ago; a more rigorous and modern version of his proof is given here. 

 Along the way, we deduce some miscellaneous results for the quadruples
\[ (2,3,2,6), \, (2,4,2,8), \, (3,2,3,9). \] 
For instance, the result for $(2,3,2,6)$ says the following: 
let $F_1, F_2$ be two general ternary cubics and $E$ a smooth planar 
cubic curve apolar to $F_1, F_2$ (in the sense explained below). Then 
$E$ passes through exactly three sextuples in $\AS$. 

In each of the cases above, there is a specific feature 
of the free resolution of $s$ general points in $\P^n$ which is 
exploited to deduce the answer. For the arguments to work smoothly, we require a 
technical condition on the polar $s$-hedra, namely that 
they be `resolution-general' (in the sense of Definition \ref{defn.resgen}).
Although the specific technique used depends on the case at hand, two general 
themes are identifiable: the geometry of associated points and 
intersection theory on symmetric products of elliptic curves. 
I do not know of any technique which would apply uniformly 
to all $(n,d,r,s)$. 

This subject is broadly referred to as `reduction to canonical form' or 
`Waring's problem for algebraic forms'; see \cite{ego.car,Dolgachev2,Ger1,IK} 
for an introduction and further references. The paper \cite{RanestadSch} is 
an excellent compendium of known results about the structure of 
$\AS$ when $r=1$. 
For a discussion of ternary cubics (the case $n=2,d=3$), 
see \cite{FranzLondon,ReichBZ}. 

\smallskip 

\noindent {\sc Acknowledgements:}
I would like to thank Anthony V.~Geramita, Leslie Roberts and 
Queen's University for financial support while this work was in progress. 
The program Macaulay-2 has been useful, and I am grateful to its 
authors Dan Grayson and Mike Stillman. 

\section{Preliminaries} \label{section.prelim}
In this section we establish notation and describe the basic set-up 
of apolarity.  
The proofs may be found in \cite{IK}, also see 
\cite{DolgachevKanev, EhRo, Ger1, Iarro3}.

The base field is $\complex$.  
Let $V$ be an $(n+1)$-dimensional $\complex$-vector space and consider the  
symmetric algebras  
\[ 
   R = \bigoplus_{i \ge 0} \, \Sym^i \, V^* , \quad   
   S = \bigoplus_{j \ge 0} \, \Sym^j \, V. 
\]   
If 
$\uu = \{u_0,\dots,u_n\}, \quad \ux = \{ x_0,\dots,x_n\}$, 
are dual bases of $V^*$ and $V$ respectively, then  
\[ R = \complex \, [u_0,\dots,u_n], \quad  
   S = \complex \, [x_0,\dots,x_n].  
\]  
There are internal product maps 
$R_i \otimes S_j \stackrel{f_{ij}}{\lra} S_{j-i}$ (see 
e.g.~\cite[p.~476]{FH}), so $S$ acquires the structure of a 
graded $R$-module.  
With the identification $u_\ell = \frac{\partial}{\partial x_\ell}$, the  
internal product can be seen as partial differentiation: 
if $\varphi \in R_i$ and $F \in S_j$, then 
$f_{ij}\, (\varphi \otimes F)$ is obtained by applying the 
differential operator 
$\varphi(\frac{\partial}{\partial x_0}, \dots , 
\frac{\partial}{\partial x_n})$ to $F(x_0, \dots, x_n)$. 
We will write $\varphi \circ F$ for $f_{ij}\, (\varphi \otimes F)$. 

Let $\Lambda \subseteq S_d$ be an $r$-dimensional subspace of degree  
$d$ forms in the $\ux$, defining a point in the Grassmannian 
$G(r, S_d)$. Let 
\begin{equation}
\Lambda^\perp = \{ \varphi \in R: \varphi \circ F = 0 \; 
\text{for every $F$ in $\Lambda$}\}. 
\end{equation}
Then $\Lambda^\perp = \bigoplus\limits_i \Lambda^\perp_i$ is  
a graded ideal in $R$, with $\Lambda^\perp_i = R_i$ for $i > d$. (It follows 
that the quotient $R/\Lambda^\perp$ is an artin level algebra of socle 
degree $d$ and type $r$, but we will not use this explicitly.) 

For $i \le d$, the codimension of $\Lambda^\perp_i$ in $R_i$ equals the 
dimension of the image of the internal product map 
\[ R_{d-i} \otimes \Lambda \lra S_i \] Hence 
\begin{equation}
\dim \Lambda^\perp_i \ge \max \, \{ 0, \dim R_i - r. \dim R_{d-i} \}. 
\label{estimate1} \end{equation}
Equality always holds for $i=d$, and it holds for all $i < d$ if 
$\Lambda$ is a general point in $G(r, S_d)$. 

We will commonly use geometric language in the sequel, e.g., 
if $n=3$, then a point in $G(2,S_4)$ will be called a pencil of 
planar quartics. 

\begin{Remark} \rm 
If $\varphi \circ F = 0$, then $\varphi, F$ were classically 
said to be apolar to each other; and sometimes the entire set-up is 
called apolarity. Of course, all of the above is subsumed in the 
statement that $R, S$ are dual Hopf algebras such that all 
structure maps are $SL(V)$-equivariant. 
\end{Remark} 

Henceforth we set $\P^n = \P S_1 = \proj R$.  
Usually $Z \subseteq \P^n$ will denote a closed subscheme with 
(saturated) ideal $I_Z \subseteq R$. 
\begin{Definition} \sl 
(cf.~\cite[Definition 5.1]{IK}) \sl 
The scheme $Z$ is said to be apolar to $\Lambda$, if  
$I_Z \subseteq \Lambda^\perp$.  
\end{Definition}  
\noindent The point of the definition is the following: 
\begin{Theorem}[Reye] \sl 
If $Z$ consists of $s$ distinct 
points $\{ L_1, \dots, L_s \} \subseteq \P^n$, then 
$Z$ is apolar to $\Lambda$ if and only if 
$\Lambda \subseteq \vspan \, \{L_1^d, \dots, L_s^d \}$.
\end{Theorem} 
We would like to consider the family of such $Z$, but for technical 
reasons, we single out those schemes whose ideals are well-behaved. 
\begin{Definition} \sl 
A (zero-dimensional) length $s$ scheme 
$Z \subseteq \P^n$ will be called resolution-general,  
if the graded Betti numbers in the minimal resolution of $I_Z$ are the 
same as those in the resolution of $s$ general points. 
\label{defn.resgen} \end{Definition} 
For instance, a length $7$ subscheme $Z \subseteq \P^2$ is 
resolution-general iff its minimal resolution looks like 
\[ 0 \ra R(-5) \oplus R(-4) \ra R(-3)^3 \ra R \ra R/I_Z \ra 0. \] 
In particular, $Z$ does not lie on a conic. 

\begin{Definition} \sl 
A zero-dimensional scheme $Z \subseteq \P^n$ will be 
called a polar polyhedron of $\Lambda$, if it is apolar to $\Lambda$ and 
resolution-general. 
\label{defn.polyhedron} \end{Definition} 
Let $\Hilb(s,\P^n)$ be the Hilbert scheme parametrising 
length $s$ subschemes of $\P^n$. Let 
$\AS(s, \Lambda)$ denote the set of polar $s$-hedra of $\Lambda$, it 
is then a constructible subset of $\Hilb(s,\P^n)$. We will write 
$\AS$ for $\AS(s, \Lambda)$ if no confusion is likely. 

\begin{Remark} \rm 
In the literature there is no unanimity on the definition of a 
`polar polyhedron', in particular the approaches in 
\cite{DolgachevKanev} and \cite{RanestadSch} are different from ours and 
from each other. It is understood that 
if $Z = \{L_1, \dots, L_s\}$ are $s$ general points, 
then morally $Z$ should count as a polar $s$-hedron of 
any $\Lambda \subseteq \vspan \, \{L_i^d\}$. However, it is not 
obvious which degenerations of $Z$ should be allowed, and it seems 
that (within reason) we should tailor our definition to the specific problem 
at hand. Many of our results depend on a free resolution of $I_Z$, and 
hence  `resolution-general' seems to be the most suitable notion. 
This issue never arises in \cite{FranzLondon}, because 
there it is tacitly assumed that all geometric 
configurations are nondegenerate. 
\end{Remark} 

If $\AS(s, \Lambda)$ is nonempty, so is $\AS(t, \Lambda)$ for 
any $t > s$. It is the case that every 
$\Lambda$ in $G(r,S_d)$ admits a polar 
$\binom{n+d}{d}$-hedron. An elementary parameter count 
(see \cite{ego.car}) shows that a general $\Lambda$ in 
$G(r, S_d)$ will admit a polar $s$-hedron only if  
\begin{equation}
s \ge \frac{r \, \binom{n+d}{d}}{n+r}. \label{formula1}  
\end{equation}

\begin{Definition} \sl 
A quadruple $(n,d,r,s)$ which satisfies  
(\ref{formula1}) is said to be nondegenerate, if a general  
$\Lambda$ admits a polar $s$-hedron. 
\end{Definition}  
A quadruple satisfying (\ref{formula1}) is degenerate if the 
set $\{\Lambda: \AS(s, \Lambda) \neq \emptyset\}$ fails to be dense 
in $G(r, S_d)$. Very few such examples are known 
(see \cite{ego.car} for the list), but none of them is without 
its geometric peculiarity. In general it is not trivial to 
prove that a particular quadruple is nondegenerate. 

For $r=1$, we have the following classification theorem by 
Alexander and Hirschowitz. 
\begin{Theorem}[see \cite{Iarro3}] \sl 
Assuming $r =1 $ and $d > 2$, the only degenerate cases are  
$(n,d,s) = (2,4,5),(3,4,9),(4,4,14)$ and $(4,3,7)$. 
\end{Theorem}
For $r> 1$ we have the following results by Dionisi 
and Fontanari. 
\begin{Theorem} \sl Assume $r > 1$. Then 
\begin{enumerate} 
\item[(i)] 
for $n=2$, the only degenerate quadruple is $(2,3,2,5)$;
\item[(ii)]
there are no degenerate quadruples with $r \ge n+1$. 
\end{enumerate} \label{theorem.fontanari} \end{Theorem} 
The proofs may be found in \cite{Fontanari2, Fontanari1} respectively. 
Part (i) was claimed by Terracini \cite{Terra1}, but his proof is obscure. 

If $(n,d,r,s)$ is nondegenerate, then with a slight  
abuse of notation we will write $\AS$ for $\AS(s, \Lambda)$, where  
$\Lambda$ is understood to be a general point of $G(r, S_d)$. It has 
dimension $s(n+r) - r \, \binom{n+d}{d}$.

\section{Associated systems of points}  
Recall (\cite[p.~313]{EGH}) that if $\Gamma$ is a 
zero-dimensional Gorenstein scheme, 
then any closed subscheme $\Gamma' \subseteq \Gamma$ has 
a residual scheme $\Gamma'' \subseteq \Gamma$, such that 
\[ \deg \Gamma' + \deg \Gamma'' = \deg \Gamma. \] 
In particular this applies if $\Gamma$ is a (global) 
complete intersection in $\P^n$, which is the only case we will need. 

Now let $\Lambda$ denote a general pencil of planar quartics.
Then $\AS = \AS(8, \Lambda)$ is $2$-dimensional; we will show 
that it is rational.  
Every $Z \in \AS$ has a Hilbert-Burch resolution 
\[ 0 \ra R(-5)^2 \stackrel{\mu}{\ra} R(-4) \oplus R(-3)^2 \ra 
R \ra R/I_Z \ra 0. \] 
(See \cite{CGO} for the basic theory behind the Hilbert-Burch theorem.) 
In particular $\dim \, (I_Z)_3 = 2$, so $Z$ has an associated point 
$\alpha(Z)$, defined to be the residual intersection of 
cubics passing through $Z$. 
The matrix of the map $\mu$ has the form 
\begin{equation} M = \left[ \begin{array}{ccc} 
{\underline 2} & {\underline 2} & {\underline 1} \\ 
{\underline 2} & {\underline 2} & {\underline 1} \end{array} \right], 
\label{HBZ} \end{equation}
with the convention that ${\underline j}$ stands for a 
degree $j$ form. 

\begin{Theorem}  \sl 
Let $\Lambda$ be a general pencil of planar quartics. Then 
the morphism $\alpha: \AS \lra \P^2$
admits a rational inverse, hence $\AS$ is a rational surface. 
\end{Theorem}  
\demo Fix a general point in the image of $\alpha$, by 
change of coordinates we assume it to be $P = [0,0,1]$. 
We would like to show that there is a \emph{unique} resolution general 
length $8$ scheme $Z$ with associated point $P$. 

Now $P$ is defined by the vanishing of the rightmost 
column in (\ref{HBZ}), hence, after row-operations, $M$ can be 
brought into the form 
\[ M = \left[ \begin{array}{ccc} 
q_1 & q_2 & u_0 \\ 
q_3 & q_4 & u_1 \end{array} \right], \quad q_i \in R_2. 
\label{HB} \]

\noindent We start with the $24$-dimensional 
vector space of $2 \times 2$ matrices 
\[ V_1 = \{ N = \left[ \begin{array}{cc} 
q_1 & q_2 \\ q_3 & q_4 \end{array}\right] : q_i \in R_2\}. \] 
For $N \in V_1$, write 
\begin{equation} 
 \theta_N = u_1 q_1 - u_0 q_3, \; \theta'_N = u_1 q_2 - u_0 q_4, \; 
   \omega_N = q_1 q_4 - q_2 q_3, 
\label{theta.omega} \end{equation}
and let $J_N$ be the ideal generated by $\theta_N, \theta'_N, \omega_N$. 
Thus $V_1$ is a parameter space for all Hilbert-Burch matrices as 
above. For a dense open set of elements $N$ in $V_1$, the ideal 
$J_N$ defines a planar length $8$ scheme. 

We let $\GL_2(\complex)$ act on $V_1$ by right multiplication, i.e., for 
$g = \left[ \begin{array}{cc} 
\alpha & \beta \\ \gamma & \delta \end{array} \right] \in \GL_2$,
and $N$ as above, 
\begin{equation} Ng = \left[ \begin{array}{cc} 
q_1 \alpha + q_2 \gamma & q_1 \beta + q_2 \delta \\ 
q_3 \alpha + q_4 \gamma & q_3 \beta + q_4 \delta \end{array} \right] 
\label{Ng} \end{equation}

Define $V_2 = \{ N \in V_1: \theta_N, \theta_N' \in \Lambda^\perp_3 \}$, 
which is a $12$-dimensional subspace of $V_1$. 
(If $F \in \Lambda$, then $\theta_N \circ F = \theta_N' \circ F = 0$ 
is a set of six linear equations. In all, $V_2$ is defined by 
$12$ linear equations which are independent for a general $\Lambda$,
hence $\dim V_2 = 12$.) 
Inside $V_2$, there is a $6$-dimensional subspace 
\[ V_3 = \{ \left[ \begin{array}{cc} 
a \, u_0 & b \, u_0 \\ a \, u_1 & b \, u_1 \end{array} \right]: 
a,b \in R_1 \}. \] 
(Since $\theta_N, \theta_N' = 0$ for $N \in V_3$, the containment 
$V_3 \subseteq V_2$ is clear.) Form the $6$-dimensional space $W = V_2/V_3$. 
For $N \in V_2$, write $[N]$ for the corresponding point in 
the projective space $\P W \simeq \P^5$. 
Since $V_3 \subseteq V_2 \subseteq V_1$ are inclusions of 
$GL_2$-modules, $W$ is also a (right) $GL_2$-module; in particular 
$PGL_2$ acts on $\P W$. 
The point of this construction lies in the following lemma: 
\begin{Lemma} \sl 
\begin{enumerate} 
\item[(i)] If $N, \tN \in V_2$ are such that $[N],[\tN]$ lie in the 
same $PGL_2$-orbit of $\P W$, then 
$J_N = J_\tN$. 
\item[(ii)] Let $Z \in \alpha^{-1}(P)$. Consider two minimal 
resolutions of $I_Z$ with corresponding Hilbert-Burch matrices 
$M, \tM$, and let $N,\tN$ denote their leftmost minors. Then 
$[N],[\tN]$ lie in the same $PGL_2$-orbit in $\P W$. 
\end{enumerate} 
\end{Lemma} 
\demo 
By straightforward calculation,
\begin{equation}
  \theta_{Ng} = \alpha \, \theta_N + \gamma \, \theta_N', \quad 
   \theta_{Ng}' = \beta \, \theta_N + \delta \, \theta_N',  \quad 
   \omega_{Ng} = \det(g) \, \omega_{N}, 
\label{ng} \end{equation}
so $J_N = J_{Ng}$. Let $Q = 
\left[ \begin{array}{cc} 
a \, u_0 & b \, u_0 \\ a \, u_1 & b \, u_1 \end{array} \right] \in V_3$. 
Then 
\begin{equation} \theta_{N +Q} = \theta_{N}, \quad 
   \theta_{N+Q}' = \theta_{N}', \quad 
   \omega_{N +Q} = \omega_N - a \, \theta_N' + b \, \theta_N, 
\label{nq} \end{equation}
so $J_{N+Q} = J_N$. This proves (i). 

Any two minimal resolutions of $I_Z$ are isomorphic 
(see \cite[\S 20.1]{Ei}), which translates into 
the statement that $N$ and some $GL_2$-translate of $\tN$ must differ 
by an element of $V_3$. This says 
that $[N],[\tN]$ must be in the same orbit, which is (ii). \qed 

Now define a subvariety 
$ Y = \{ [N] \in \P W: \omega_N \circ \Lambda = 0\} \subseteq \P W$.
Formulae (\ref{nq}) imply that 
$\omega_{N+Q} \circ \Lambda = 0 \iff \omega_N \circ \Lambda = 0$
(since $\theta_N \circ \Lambda = \theta_N' \circ \Lambda = 0$), 
hence this definition is meaningful. The inclusion 
$Y \subseteq \P W$ is a $PGL_2$-stable by formulae (\ref{ng}). By the previous 
lemma, each $Z \in \alpha^{-1}(P)$ defines an orbit 
$\Omega_Z \subseteq Y$.  The $PGL_2$-stabilizer of a point in 
$\Omega_Z$ is trivial, hence $\dim \Omega_Z = 3$. The union 
of $\{\Omega_Z\}_{Z \in \alpha^{-1}(P)}$ fills a dense open subset 
in $Y$. Hence it is enough to show that $Y$ contains only one 
three-dimensional component, this will imply that 
$\alpha^{-1}(P)$ is singleton. Define 
\[ \begin{aligned} 
\Gamma_1 & = \{ [N] \in \P W: N = \left[ \begin{array}{cc} 
q_1 & 0 \\ q_3 & 0 \end{array} \right] \text{for some $q_i$ and 
$\theta_N \circ \Lambda = 0$} \}, \\ 
\Gamma_2 & = \{ [N] \in \P W: N = \left[ \begin{array}{cc} 
0 & q_2 \\ 0 & q_4 \end{array} \right] \text{for some $q_i$ and 
$\theta_N' \circ \Lambda = 0$} \},
\end{aligned} \] 
each of which is a copy of $\P^2$ in $Y$. Define a birational map 
$h: \Gamma_1 \lra \Gamma_2$ as follows. Let $[N] \in \Gamma_1$, then 
there is a $4$-dimensional 
family of solutions $(q_2,q_4)$ to the equations 
\[ \theta_N' \circ \Lambda = \omega_N \circ \Lambda = 0. \] 
(This is so because $q_2,q_4$ together depend upon $12$ parameters and 
there are $8$ equations.) However, if $(q_2,q_4)$ is one such 
solution, then $(q_2+ au_0,q_4+ au_1)$ is also one for any 
$a \in R_1$, and this accounts for all the solutions. Hence the 
class in $\P W$ of the matrix $\left[ \begin{array}{cc} 
0 & q_2 \\ 0 & q_4 \end{array} \right]$ is uniquely determined. 
We define $h([N])$ to be this class. (The reader should verify that 
this definition is independent of the choice of coset representative 
for $[N]$.) Now a general element in $Y$ can be written as a sum 
$[N] + [h(N)]$ for $[N] \in \Gamma_1$, i.e., 
the ruled join of $\Gamma_1, \Gamma_2$ along $h$ contains 
a dense open subset of $Y$. Since this join is irreducible (it is the 
image of the Segre imbedding $\P^2 \times \P^1 \subseteq \P^5$), 
we are done. \qed 

The argument for the following proposition is similar. 
As before, $(2,3,4,7)$ is nondegenerate by Theorem 
\ref{theorem.fontanari}. 
\begin{Proposition} \sl 
Let $\Lambda$ be a general web of planar cubics. Then 
$\AS(7, \Lambda)$ is a rational surface. 
\end{Proposition} 
\demo 
The Hilbert-Burch matrix for $Z \in \AS$ is 
$\left[ \begin{array}{ccc} 
\underline{1} & \underline{1} & \underline{1} \\ 
\underline{2} & \underline{2} & \underline{2} \end{array}\right]$. 
For a general $Z$, the linear forms in the top row are independent, 
hence after column operations we can assume the matrix to be 
\[ \left[ \begin{array}{ccc} 
u_0 & u_1 & u_2 \\ 
q_0 & q_1 & q_2 \end{array}\right], \quad q_i \in R_2. \] 
Let $V_1$ denote the $18$-dimensional vector space
$\{ [q_0,q_1,q_2] : q_i \in R_2\}$, 
and $V_2$ the $3$-dimensional subspace 
$\{ [ a u_0, a u_1,a u_2] : a \in R_1\}$. 
Let $W = V_1/V_2$. Then the $12$ equations 
$\{ (u_i q_j - u_j q_i) \circ \Lambda = 0\}$ cut out a $2$-plane 
in $\P W$ which is birational to $\AS$. \qed 

\smallskip 

Now let $(n,d,r,s) = (3,2,6,7)$, we will show 
that $\AS$ is birational to the projective $3$-space.  
The ideal of every $Z \in \AS$ is generated by three quadrics and 
and a cubic. The associated point $\alpha(Z)$ is defined 
to the residual intersection of the quadrics through $Z$.  

\begin{Proposition} \sl 
Let $\Lambda$ be a general point of $G(6, S_2)$. 
Then the map $\alpha: \AS \lra \P^3$ is birational. 
\end{Proposition}  
\demo  
Let $Z$ be a resolution-general scheme of length $7$. 
It is apolar to $\Lambda$ iff the three generating 
quadrics lie in $\Lambda^\perp_2$. 

Let $P$ be a general point of $\P^3$, and let 
$W \subseteq \Lambda^\perp_2$ be the $3$-dimensional subspace of 
forms vanishing at $P$. Then $W$ defines a length $8$ scheme $Y$. 
Now the residual scheme of $P$ in $Y$ is the 
only point of $\AS$ mapping to $P$. \qed 

\begin{Remark} \rm  
The case $(2,3,7,8)$ has a similar geometry, where $\AS$ is 
birational to $\P^2$. For $(2,3,8,8)$ (resp.~$(3,2,7,7)$), $\AS$ is 
a finite set consisting of $9$ (resp.~$8$) points. 
\end{Remark}  

\section{Symmetric Products of Elliptic Curves}
For the examples in this section, 
the determination of $\AS$ reduces to an intersection-theoretic 
calculation on the symmetric product of an elliptic curve. 
If $E$ is a smooth projective curve, then $E^{(m)}$ will denote  
its $m$-th symmetric product. This is a smooth projective variety 
whose points are naturally seen as effective degree $m$ divisors on $E$.  

Let $\Lambda$ be a general net of planar quartics. Since 
$(2,4,3,9)$ is nondegenerate, $\AS$ is a finite set. In the next 
theorem we calculate its cardinality.  

\begin{Theorem}  \sl 
Let $\Lambda$ be a general net of planar quartics. Then 
$\Lambda$ admits $4$ polar enneahedra. 
\label{prop.count.2439} \end{Theorem}  
\demo  
The ideal of $Z \in \AS$ is generated by one 
cubic and $3$ quartics. 
The space $\Lambda^\perp_3$ is one-dimensional, i.e., $\Lambda$ 
is apolar to a unique cubic curve $E \subseteq \P^2$. 
Since $\Lambda$ is general, we may (and will) assume 
that $E$ is smooth. If $H$ denotes the hyperplane divisor 
on $E$, then we have an identification 
$H^0(E, 4H) =  R_4/(I_E)_4$. This is a $12$-dimensional 
space, denoted $U$.  

Let $W = \Lambda^\perp_4/(I_E)_4$, which is a $9$-dimensional 
space inside $U$. 
Every scheme $Z \subseteq \P^2$ of length $9$ which is apolar to  
$\Lambda$ is contained in $E$, and thus defines an effective  
divisor on $E$. Then the $3$-dimensional space 
$H^0(E, 4H - Z)$, which is \emph{a priori} inside  
$U$, is in fact contained in $W$. The argument shows that  
the following diagram is a fibre square: 
\[ \begin{CD}
    \AS       @>>>  G(3,W) \\  
    @VVV      @V{i_1}VV    \\
    E^{(9)}  @>{i_2}>> G(3,U) \end{CD} 
\]   
Here $i_1$ is the natural inclusion and 
$i_2(Z) = H^0(E,4H - Z)$. Since the images of both inclusions 
have complementary codimensions,  it is enough to take 
the intersection of their classes inside $H^*(G(3, U), \Z)$ in 
order to calculate the degree of $\AS$ as a zero-cycle.  

\medskip 

\noindent {\bf Conventions.} 
The notation for Schubert calculus follows \cite[\S 14.7]{Fu1}.
We refer to \cite{ACGH} for some basic cohomological calculations 
on curves. If $X_1, X_2$ are varieties, then denote projections 
by $\pi_i: X_1 \times X_2 \lra X_i$. 
All cohomology is with $\Z$-coefficients. If $\alpha$ is a class in 
$H^*(X_1)$ (resp.~$H^*(X_2)$), then its pullback to $H^*(X_1 \times X_2)$ is 
denoted $\alpha \otimes 1$ (resp.~$1 \otimes \alpha$). Cup product 
is written as juxtaposition. 

Firstly, we should find the rank $3$ subbundle of $U \otimes \O_{E^{(9)}}$ which defines 
the inclusion $i_2$. Let $\Delta$ denote the  
universal divisor on $E^{(9)} \times E$ (see \cite[Ch. IV]{ACGH}), so that  
$\Delta|_{\{Z\} \times E} = Z \times E$. Define a line bundle 
$\MM = \pi_2^*(\O_E(4H)) \otimes \O(-\Delta)$  on $E^{(9)} \times E$. Applying 
${\pi_1}_*$ to the inclusion 
\[ \MM \subseteq \pi_2^*(\O_E(4H)), \] 
we have 
\[ (\GG =) {\pi_1}_* (\MM) \subseteq  U \otimes \O_{E^{(9)}}. \] 
A moment's reflection will show that 
$i_2$ is induced by the last inclusion. 

The image of $i_2$ has class $\{3,3,3\}$. Hence by 
the Jacobi-Trudi identity, the class of $\AS$ in $H^{18}(E^{(9)})$ is 
given by $c_3(\GG^*)^3$, which we now calculate. 
\smallskip 

\noindent {\bf The cohomology rings of $E$ and $E^{(9)}$.} 
Let $\delta_1, \delta_2 \in H^1(E)$ be 
a symplectic basis, it will then generate $H^*(E)$. The product 
$\eta = \delta_1 \, \delta_2 \in H^2(E)$ is the class of a point. 

Let $\LL$ be a Poincar{\'e} line bundle 
(\cite[Ch. IV]{ACGH}) on $E \times \Pic^9(E)$, then $\EE = {\pi_2}_*(\LL)$ 
is a rank $9$ bundle on $\Pic^9(E)$. Fix an isomorphism $\Pic^9(E) = E$, then 
by the calculation of \cite[p.~336]{ACGH}, $c_1(\EE) = - \eta$. Now 
let $\xi = c_1(\O_{\P \EE}(1)) \in H^2(\P \EE)$. With the identification 
$\P \EE = E^{(9)}$, 
the ring $H^*(E^{(9)})$ is generated by $\xi$ and (the pullbacks of) 
$\delta_1, \delta_2$, subject to the relation $\xi^9 = \xi^8 \eta$. 

\noindent{\bf The Chern class of $\MM$ and G-R-R.} Let 
\[ - \gamma = (\delta_1 \otimes 1)(1 \otimes \delta_2) 
          - (\delta_2 \otimes 1)(1 \otimes \delta_1),  \] 
a class in $H^{1,1}(E^{(9)} \times E)$. By \cite[p.~337-338]{ACGH}, 
\[ c_1(\O(\Delta)) = \xi \otimes 1 + \gamma + 9 (1 \otimes \eta), \] 
hence 
\[ c_1(\MM) =  - \, \xi \otimes 1 -  \gamma + 3 (1 \otimes \eta). \] 
Now we apply Grothendieck--Riemann--Roch to $\MM$ along the projection 
$E^{(9)} \times E \stackrel{{\pi_1}}{\lra} E^{(9)}$. Thus 
\[ 
\ch({\pi_1}_! \, \MM)\, \td(E^{(9)}) = 
{\pi_1}_* (\ch(\MM) \, \td(E^{(9)} \times E)). 
\]
Since $R^i{\pi_1}_* \MM = 0$ for $i > 0$ and $\td(E) = 1$, this simplifies to 
\[ 
\ch (\GG) = {\pi_1}_*(e^{c_1(\MM)}). \label{formula.GRR} 
\] 
Let $n_i$ denote the $i$-th Newton class of $\GG$ 
(i.e., the sum of $i$-th powers of the Chern roots of $\GG$), 
then $\ch(\GG) = \sum\limits_{i \ge 0} {n_i}/{i!}$. 
Now we expand the exponential series, and apply ${\pi_1}_*$ 
term by term, to get 
\[ \begin{aligned} 
n_0 = 3,  \quad & n_1 = \frac{1}{2}(-6 \, \xi - 2 \eta), \\ 
n_2 = \frac{1}{3} (9 \, \xi^2 + 6 \, \xi \eta), \quad & n_3 = 
\frac{1}{4}(-12 \, \xi^3 - 12 \, \xi^2 \eta). 
\end{aligned} \] 
Then \[ 
c_3(\GG) = \frac{1}{6}n_1^3 - \frac{1}{2} n_1 n_2 + \frac{1}{3} n_3 
= - (\xi^3 +\xi ^2 \eta). \] 
Hence finally 
\[ c_3(\GG^*)^3 = (\xi^3 + \xi ^2 \eta)^3 = 4 \, \xi^8 \eta. \] 
Since $\xi^8 \eta$ is the class of a point on $E^{(9)}$, we deduce that 
$\AS$ has degree $4$. 

In order to show that $\AS$ is reduced and hence consists of $4$ geometric 
points, we use Kleiman's transversality result (see \cite[Theorem 10.8]{Ha}). 
We can reformulate the entire construction in the following way: start 
with a smooth $E$ and hence $U$, then specifying a 
codimension $3$ subspace $W \subseteq U$ is tantamount to 
specifying $\Lambda$. Since $G(3,U)$ is a homogeneous space for $GL(U)$, 
the intersection is transversal for a general $W$, so $\AS$ is reduced. 
\qed 

\smallskip 

The next example is that of a pencil of cubic surfaces. 
We need to show that $(3,3,2,8)$ is nondegenerate, the 
proof is given in \S \ref{grove.3328}. 

\begin{Proposition} \sl 
Let $\Lambda$ be a general pencil of cubic surfaces. Then 
$\Lambda$ admits $3$ polar octahedra. 
\end{Proposition}  
\demo  
The calculation is very similar to Theorem \ref{prop.count.2439}.  
The ideal of $8$ general points in $\P^3$ is generated by $2$ quadrics 
and $4$ cubics. 
Now $\Lambda^\perp_2$ is $2$-dimensional, hence generates the 
ideal of a smooth normal elliptic quartic $ E \subseteq \P^3$ 
apolar to $\Lambda$, and every $Z \in \AS$ is in fact 
contained in $E$. Let 
\[ U = R_3/(I_E)_3,  \quad   
W = \Lambda^\perp_3/(I_E)_3, \] which are 
spaces of dimension $12, 10$ respectively. Define $i_1, i_2$ as before, then 
the following diagram is a fibre square 
\[ \begin{CD}
    \AS       @>>>  G(4,W) \\  
    @VVV      @V{i_1}VV    \\
    E^{(8)}  @>{i_2}>> G(4,U) \end{CD} 
\]   
Now $i_2$ is induced by a rank $4$ bundle $\GG$ on $E^{(8)}$. 
The class of $\AS$ in $E^{(8)}$ equals 
\[ c_4(\GG^*)^2 = (\xi^4 + \xi^3 \eta)^2 = 3 \, \xi^7 \eta. \] 
The argument for transversality is the same as before. 
\qed  

\smallskip 

Using similar calculations, 
we can give alternate proofs of the following results 
by Schlesinger \cite[p.~212]{Sch1}). The original 
argument uses $\vartheta$-functions.  
\begin{Proposition}[Schlesinger] \sl 
\begin{enumerate} 
\item 
Let $\Lambda$ be a general pencil of planar cubics.  Fix a 
general elliptic curve $E \subseteq \P^2$ apolar to $\Lambda$. Then 
there are $3$ polar hexahedra of $\Lambda$ which are 
contained in $E$. 
\item 
Let $\Lambda$ be a general pencil of planar quartics. Fix 
a general elliptic curve $E \subseteq \P^2$ apolar 
to $\Lambda$. Then there are $3$ polar octahedra of $\Lambda$ 
which are contained in $E$. 
\end{enumerate} \end{Proposition}  
\demo  
We will only prove (1), the argument for (2) is identical in essence. 
Recall that the ideal of $6$ general planar points is generated 
by $4$ cubics. 
Since $(2,3,2,6)$ is nondegenerate\footnote{This is the smallest 
$s$ possible, because $(2,3,2,5)$ is degenerate by \cite{ego.car}.}, 
$\AS(6, \Lambda)$ is $4$-dimensional. 
Consider the incidence correspondence 
\[ \Phi \subseteq \AS \times \P \Lambda^\perp_3, 
\quad \Phi = \{ (Z, E): Z \subseteq E \}. \] 
The projection $\pi_1: \Phi \lra \AS$ is generically a $\P^3$-bundle, so 
$\dim \Phi = 7$. 
Fix a general elliptic curve $E$ apolar to $\Lambda$, and consider the diagram 
\[ \begin{CD}
    {}     @.            G(3,\Lambda^\perp_3/(I_E)_3) \\  
    @.      @V{i_1}VV    \\
    E^{(6)}  @>{i_2}>> G(3,R_3/(I_E)_3) 
\end{CD} \]   
As usual, $i_1$ is the inclusion and $i_2(Z) = H^0(E,3H - Z)$. Then 
$i_2(Z)$ lies in the image of $i_1$, iff $Z$ is apolar to $\Lambda$. 
Calculating as before, the product 
$[\image \, i_1].[\image \, i_2]$ equals thrice the class of a point. 
Hence $\pi_2^{-1}(E)$ must be nonempty. 
This implies that 
$\pi_2: \Phi \lra \P \Lambda_3^\perp (\simeq \P^7)$ is dominant. 
But then it is generically finite, hence for a general $E$, the fibre 
$\pi_2^{-1}(E)$ consists of $3$ points. 
\qed 

\smallskip 

It is shown in \cite{ego.car} (using a machine calculation) that 
$(5,2,3,9)$ is nondegenerate. Now there is a (unique) 
elliptic sextic curve passing through $9$ general points of 
$\P^5$. (The classical reference is \cite{Room1}, also see~\cite{Dolgachev1} 
for a proof using Gale duality.) Hence 
if $\Lambda$ is a general net of quadrics in $\P^5$ and $Z$ a set of 
$9$ general points apolar to $\Lambda$, then the elliptic 
sextic passing through $Z$ is apolar to $\Lambda$. 

\begin{Proposition} \sl 
Let $\Lambda$ be a general net of quadrics in $\P^5$. Fix a 
general elliptic sextic curve $E \subseteq \P^5$ apolar 
to $\Lambda$. Then there are $4$ polar enneahedra 
of $\Lambda$ which are contained in $E$. 
\end{Proposition} 
\demo Similar to above. Use the fact that the ideal of 
$9$ general points (resp.~an elliptic sextic curve) is generated by 
$12$ (resp.~$9$) quadrics. \qed 

\section{The $(2,3,3,6)$ case} 
Now we come to London's beautiful calculation in \cite{FranzLondon}, where 
he determines the number of polar hexahedra of a general net 
of cubic curves. I have rewritten the proof so as to make it more 
transparent, but all the key ideas are already in the original. 

Let $\Lambda$ be such a net. By Theorem \ref{theorem.fontanari}(i), 
$\Lambda$ has a finite number of polar hexahedra. We will count 
them by setting up a correspondence on a certain elliptic curve. 

\subsection{} \label{subsection.motivation}
We begin by motivating the constructions which are to follow. 
Say $\{F_1, F_2, F_3\}$ is a basis of $\Lambda$ and $Z = 
\{L_1, \dots, L_6\}$ one of its polar hexahedra. 
We have expressions 
\[ F_j = c_{1j} \, L_1^3 + \dots + c_{6j} \, L_6^3, \quad j=1,2,3. \] 
Let $\psi \in R_2$ be the form which defines the conic passing 
through $\{L_2, \dots, L_6\} \subseteq \P S_1$. Since 
$\psi$ annihilates $L_2^3, \dots,L_6^3$, we have 
$\psi \circ F_j = \text{constant} \times L_1$ for every $j$, 
so $\psi \circ \Lambda$ is only a $1$-dimensional vector space. 
It will be seen below (\S \ref{defn.Psi}) that all $\psi$ with this 
property lie on a curve. Similarly if 
$l_1,l_1' \in R_1$ annihilate $L_1$, then the six 
derivatives $\{l_1 \circ F_j, l_1' \circ F_j: j=1,2,3 \}$ span 
only a $5$-dimensional space. It will be seen below (\S \ref{defn.E}) 
that all $2$-dimensional spaces $\vspan \{l_1, l_1' \} \subseteq R_1$ 
with this property lie on a curve, isomorphic to the previous one. 

\subsection{} \label{defn.Psi}
Now we come to the actual constructions. 
The symbol (\chk ) will appear frequently, it is explained in 
Remark \ref{remark.check}. 
Consider the vector bundle morphism on $\P R_2 \, ( = \P^5)$ 
\[ f_{23}: \O_{\P^5}(-1) \otimes \Lambda \lra S_1 \] 
coming from the internal product map of \S \ref{section.prelim}. 
Define the degeneracy locus $\Psi = \{ \rank f_{23} \le 1 \}$. For 
a general $\Lambda$, it is a degree $6$ normal elliptic curve in 
$\P^5$ (\chk). Note that $\Lambda^\perp_2 = 0$ by the 
generality of $\Lambda$, so $\rank f_{23}$ is exactly $1$ at 
each $\psi \in \Psi$. 

\subsection{} \label{defn.E} 
Now identify $\P S_1$ with the Grassmannian $G(2,R_1)$, the latter 
is equipped with a rank two tautological bundle 
$\BB \subseteq R_1 \otimes \O_G$. The internal product 
$f_{13}$ gives a morphism 
\[ f_{13}: \BB \otimes \Lambda \lra S_2 \] 
The locus $E = \{ \rank f_{13} \le 5 \} = \{ \det f_{13} =0 \}$
is given by a section of $\O_{\P S_1}(3)$, hence it is a 
smooth (\chk) degree $3$ curve in $\P S_1$. By the generality 
of $\Lambda$, the rank of $f_{13}$ is exactly $5$ at 
every $L \in E$ (\chk). 

\subsection{} We have an isomorphism 
\[ \alpha: E \lra \Psi \] 
defined as follows: let $L \in E$, and $U = L^\perp$. 
By hypothesis, the space 
$f_{13} (U \otimes \Lambda)$ is $5$-dimensional, so it 
is annihilated by a unique form in $\P R_2$, we declare 
$\alpha(L)$ to be this form. It is clear 
that $f_{23}(\alpha(L) \otimes \Lambda)$ is only $1$-dimensional (since 
$U$ annihilates it), so $\alpha(L) \in \Psi$. 

If $Z$ is as in \S \ref{subsection.motivation} above, then 
$\alpha(L_1)$ is the conic envelope containing the lines 
defined by $L_2, \dots, L_6$. 

\subsection{} Define a correspondence $\T$ on $E$ as follows: 
$(L,M) \in \T$ iff $M$ lies on the conic defined by $\alpha(L)$. 
For a fixed $L$, there are $6$ positions of $M$ such that 
$(L,M) \in \T$. For a fixed $M$, the elements of $\Psi$ which 
vanish at $M$ lie on a hyperplane section of $\Psi$. Via 
$\alpha^{-1}$, the points of this hyperplane section correspond 
to $6$ positions of $L$. This shows that $\T$ has degree $(6,6)$ and 
valence zero. 

\subsection{} 
By the general theory of correspondences (see \cite[\S 2.5]{GH}), 
there are $12$ elements in $\T$ of the form 
$(L,L)$, they are called the united points of $\T$. 
Moreover $\T,\T^{-1}$ have $72$ common points, i.e., 
pairs $(L,M)$ such that $(L,M), (M,L) \in \T$. Hence there 
are $72 - 12 = 60$ such pairs where $(L,M)$ are distinct. 

It is clear that starting from $Z$, the pairs $(L_1, L_2)$ etc.~are 
common to $\T, \T^{-1}$. The next lemma says that the implication is
reversible. 
\begin{Lemma} \sl 
Assume $(L,M), (M,L) \in \T$, and $L \neq M$. Let the conics 
$\alpha(L), \alpha(M)$ intersect in $\{P_1, P_2,P_3,P_4 \}$. 
Then $Z = \{L,M,P_1,\dots, P_4 \}$ is a polar hexahedron of $\Lambda$. 
\end{Lemma} 
\demo Recall that the ideal of $6$ general points is generated 
by $4$ cubics. Let $l,l' \in R_1$ be generators of $L^\perp$, 
and $m,m'$ of $M^\perp$. Consider the four cubic forms 
\[ 
\{ l \, \alpha(L), l' \, \alpha(L), m \, \alpha(M), m' \, \alpha(M) \}.
\] 
They are linearly independent and each of them vanishes 
at all points of $Z$. Hence together they generate $(I_Z)_3$. 
Moreover, the definition of $\alpha$ implies that each of them 
annihilates $\Lambda$. Hence $I_Z \subseteq \Lambda^\perp$. \qed 

Now a polar hexahedron of $\Lambda$ gives $2 \binom{6}{2} = 30$ pairs 
$(L_i,L_j)$ common to $\T, \T^{-1}$. Alternately, starting from a 
common point we can reconstruct a polar hexahedron as shown above. 
Hence, following London, we conclude that 
$\Lambda$ has $60 \div 30 = 2$ polar hexahedra. 

\begin{Remark} \rm \label{remark.check} 
At several points in the proof we need to argue that 
our construction satisfies a certain `good' property, for instance
$\Psi$ has codimension $4$ as expected and is smooth. This 
follows from the generality of $\Lambda$, as soon as we 
verify that it holds for a specific $\Lambda$. Such points are 
marked by (\chk), and I have verified the property in question 
by a direct computer calculation for a general net in the span of 
\[ x_0^3, x_1^3, x_2^3, (x_0 + x_1 + x_2)^3, (x_0 - x_1 + x_2)^3, 
   (x_0 - 2 x_1 + 3x_2)^3. \] 
This was carried out in Macaulay-2. For instance, to verify the last point 
in \S \ref{defn.E}, we choose two basis elements with indeterminate entries 
for an element of $G(2, R_1)$, represent $f_{13}$ by a matrix and 
check that the ideal defined by all $5 \times 5$ minors defines the 
empty scheme. \end{Remark} 

\section{Nondegeneracy of $(3,3,2,8)$} \label{grove.3328}
To prove this result, we will use the notion of a {\sl grove}, which 
was introduced in \cite{ego.car}. The general definition is meaningful 
for any $(n,d,r,s)$, but we will formulate it only for the case at hand. 

Let $\up = \{p_1, \dots, p_8\}$ (resp.~$\uQ = \{Q_1, \dots, Q_8\}$) be 
points in $\P^1$ (resp.~in $\P^3$). 
\begin{Definition} \sl A grove for the data $\up, \uQ$ is 
a linear system $\Gamma \subseteq \P H^0(\P^3, \O_\P(3))$ of projective
dimension (say) $t$, satisfying the following conditions: 
\begin{enumerate} 
\item 
The base locus of $\Gamma$ contains all the $Q_i$, 
\item 
$t=0$ or $1$, and
\item 
either $t=0$ and the generator of $\Gamma$ is singular at all 
$Q_i$, or 
$t=1$ and there is an isomorphism $\gamma: \P^1 \lra \Gamma$ such 
that for every $i$, the hypersurface 
$\gamma(p_i)$ is singular at $Q_i$. 
\end{enumerate} 
\end{Definition} 

Now \cite[Theorem 2.6]{ego.car} says the following: the quadruple 
$(3,3,2,8)$ is nondegenerate iff there does not exist a 
grove for \emph{general} points $\up, \uQ$ as above. 
Existence of a grove is an open property of $\up, \uQ$ (loc.~cit.), 
so it is enough to exhibit some collection of points which does not 
admit a grove. I concede that the definition of a grove is 
awkward, in defence one can only say that it is a proof-generated 
concept in the sense of Lakatos (see \cite[Appendix 2]{Lakatos}). 
We  begin with a preliminary lemma. 
\begin{Lemma} \sl 
Let $E$ be an elliptic curve and ${\mathcal M}$ a line  
bundle on $E$ of degree $4$. Let $Q_1, \dots, Q_8$ be  
distinct points on $E$. Then it is possible to find points $p_1, \dots, p_8$  
on $\P^1$, such that there is no morphism $f: E \lra \P^1$  
satisfying the following conditions:  
\begin{itemize}  
\item[A.]
$2 \le \deg f \le 4$, and if $\deg f = 4$ then  
$f^*(\O_{\P^1}(1)) \simeq {\mathcal M}$;  
\item[B.]
the equality $f(Q_i) = p_i$ holds for at least $4 + \deg f$ values  
of $i$.  
\end{itemize} \label{lemma1} \end{Lemma}  
\demo Since $h^0({\mathcal M}) = 4$, there are $\infty^4$ $g^1_4$'s  
coming from ${\mathcal M}$. However, modulo automorphisms of $\P^1$ there are  
$\infty^5$ octuples $(p_1, \dots, p_8)$. Hence for a general octuple,  
there is no such map of degree $4$.  

Similarly there are $\infty^3$ (resp.~$\infty^1$) $g^1_3$'s (resp. $g^1_2$'s)  
on $E$. Since (B) imposes $4$ (resp.~$3$) conditions in these cases, for a  
general choice of $p_i$ none of the possibilities can hold. The 
lemma is proved. \qed  

Now $w,x,y,z$ be the coordinates in $\P^3$. Consider 
the normal elliptic quartic $E \subseteq \P^9$ defined by the two quadrics 
\[ G_1 = wx+ xy + yz + zw, \quad 
   G_2 = wy + xz. \] 
Choose points 
\[ \begin{array}{llll}
Q_1 = [1,0,0,0],  &  Q_3 = [0,0,1,0],   & Q_5 = [-1,1,1,1],  & Q_7 = [1,1,-1,1], \\ 
Q_2 = [0,1,0,0],  &  Q_4 = [0,0,0,1],   & Q_6 = [1,-1,1,1],  & Q_8 = [1,1,1,-1], 
\end{array} \]
all lying on $E$, and the $p_i = [p_{i1},p_{i2}]$ general in $\P^1$. 

Assume by way of contradiction that $\Gamma$ is a grove 
for the data. If $t=0$, then the generator of $\Gamma$ contains at 
least $16$ points of $E$ (counting each $Q_i$ as two), hence it 
contains $E$ by B{\'e}zout's theorem. 

Case 1. Assume that $\Gamma$ contains $E$ as a fixed component (with 
$t$ possibly $0$ or $1$). Then 
$\Gamma$ is spanned by two cubics of the form 
\[ C_1 = L_1 G_1 + L_2 G_2, \quad C_2 = L_1' G_1 + L_2' G_2, \] 
where $L_1,L_1'$ etc are linear forms and $p_{i1}C_1 + p_{i2}C_2$ 
is singular at $Q_i$ for $i=1, \dots, 8$. (The case $C_1 = $ (constant).$C_2$ 
corresponds to $t=0$.) An elementary linear algebra computation 
on the Jacobian matrix shows that this is impossible for general $p_i$. 
\smallskip 

Case 2. Assume that $E$ is not contained in the base locus of $\Gamma$ 
(hence necessarily $t=1$). 
Let $\lambda$ be the linear series obtained by restricting $\Gamma$ to 
$E$ and removing the base divisor $\sum Q_i$. Thus $\lambda$ is a $g^1_4$. 
Let $f : E \lra \P^1$ be the corresponding morphism (of course, 
only well-defined up to automorphisms of $\P^1$). Let $H$ denote the 
hyperplane divisor on $E$ and $\MM = \O_E(4H - \sum Q_i)$. 

Case 2.1. If $\lambda$ is base point free (i.e., if $\Gamma$ has no  
additional base point on $E$ away from $\sum Q_i$),  
then $\deg f = 4$ and $f^*(\O_{\P^1}(1)) \simeq {\mathcal M}$.  
Since the quartic $\gamma(p_i)$ passes doubly through $Q_i$,  
we have $f(Q_i) = p_i$ for all $i$.  

Case 2.2. If $\lambda$ has base points, then $\deg f \le 3$. The base locus of  
$\lambda$ can contain at most $4 - \deg f$ points from the set $\{Q_i\}$, hence  
$f(Q_i) = p_i$ holds for at least $4 + \deg f$ values of $i$.  

Now the previous lemma implies that either subcase is impossible for general 
choice of $p_i$, hence no such grove can exist. 
We have proved that $(3,3,2,8)$ is nondegenerate. 
\qed 
\section{Open problems} 
Whenever $\AS$ is a finite set, we have the obvious 
enumerative problem of counting its cardinality. Beyond a handful of 
cases (see \cite{RanestadSch}) it is entirely open. 
In particular, I do not know the cardinality of $\AS$ 
for $(2,4,4,10)$ or $(3,3,3,10)$. 

It is also of interest to consider the family of positive dimensional 
schemes (with a fixed Hilbert polynomial) apolar to $\Lambda$. For 
instance, it is known that there are two twisted cubics apolar 
to a general web of quadrics in $\P^3$ (see \cite[p.~32]{Norway2}). 

It is known that a general net of quadrics in $\P^5$ does not admit a 
polar octahedron (see \cite{ego.car}), contrary to what one would expect by 
counting parameters. However it is not known if such a net always admits an 
apolar rational normal quintic curve. A solution to this 
would help in elucidating the case $(5,2,3,8)$. 

\bibliographystyle{plain}

\vspace{1.5cm} 
\parbox{12cm}{Jaydeep V.~Chipalkatti \\ 
Department of Mathematics \\ 
University of Manitoba \\
Winnipeg, MB R3T 2N2, Canada. \\ 
email:{\tt chipalka@cc.umanitoba.ca}}

\end{document}